\documentclass{article}
\usepackage{amsmath,amssymb,amsthm, amsfonts,bm,float,inputenc}
\usepackage[pdftex]{graphicx}
\usepackage[shortlabels]{enumitem}    
 \usepackage{hyperref}

\usepackage{tikz-cd}
\usepackage{mathrsfs}     

\theoremstyle{definition}
\swapnumbers                

\newtheorem{theorem}{Theorem}
\newtheorem{definition}[theorem]{Definition}
\newtheorem{proposition}[theorem]{Proposition}
\newtheorem{lemma}[theorem]{Lemma}

\newtheorem{HRThm}[theorem]{Hurwitz-Radon Theorem}

\newcommand{\vsp}{\vspace*{6mm}}
\newcommand{\spz}{\hspace*{.5cm}}
\newcommand{\spa}{\hspace*{1.0cm}}
\newcommand{\spb}{\hspace*{2cm}}
\newcommand{\spc}{\hspace*{4cm}}
\newcommand{\spd}{\hspace*{5cm}}

\newcommand{\tr}{\textsf{T}}

\def\a{\alpha}
\def\b{\beta}

\raggedright

\begin{document}
\begin{center}
{\bf Construction of Identities for Sums of Squares}\\[2pt]
{\small Daniel B. Shapiro\\[-2pt] Ohio State University}
\end{center}

{\sc Abstract.}   From sum-of-squares formulas of sizes $[r, s, n]$ and $[r', s', n']$ we construct a formula of \\
size $[r+r',\; 2ss', \;2nn']$.

\bigskip\hrule\medskip

A sum-of-squares formula of size $[r, s, n]$ is an equation of the type:\\[5pt]
\spc  $   \left(x_1^2 + \cdots + x_r^2\right)\cdot \left(y_1^2 + \cdots + y_s^2\right) \;=\; z_1^2 + \cdots + z_n^2,$    \hfill(*)    \\[5pt]
where $X = (x_1, \dots , x_r)$ and $Y = (y_1, \dots, y_s)$ are systems of independent indeterminates, and each \linebreak 
$z_k = z_k(X, Y)$ is a bilinear form in $X$ and $Y$ with coefficients in a given field $K$.  Most constructions of such formulas use coefficients only from $\{0, 1, -1\}$. For example, the $[2, 2, 2]$-formula 
 \[   \left(x_1^2 +  x_2^2\right)\cdot \left(y_1^2 + y_2^2\right) \;=\; \left(x_1y_1\rule{0mm}{3mm}  - x_2y_2\right)^2 + \left(\rule{0mm}{3mm}x_1y_2 + x_2y_1\right)^2   \]
arises from multiplication of complex numbers.  Set $\a = x_1 + x_2i$ and $\b = y_1 + y_2i$ and note that the norm property $|\a\b| = |\a|\cdot |\b|$ yields the formula displayed above. Here $z_1 = x_1y_1  - x_2y_2$ and $z_2 = x_1y_2 + x_2y_1$ are bilinear forms in $X$ and $Y$.  

In 1748 Euler recorded a $[4, 4, 4]$-formula, and in 1843 Hamilton interpreted that as the norm property for the 4-dimensional algebra of quaternions.  A few years later, Graves and Cayley discovered the algebra of octonions and noted that its multiplication yields an $[8,8,8]$-formula.  After other mathematicians were unable to find a 16-square identity (that is, a [16, 16, 16]-formula), Adolf Hurwitz \cite{aH98} settled the question in 1898.

\theorem[Hurwitz]  If an $[n, n, n]$-formula exists, then $n = 1, 2, 4$ or 8.

In that same paper, he asked: \\[5pt]
\spa \textbf{Hurwitz Problem.} For which $r, s, n$ does there exist an $[r, s, n]$-formula?

Later, Hurwitz \cite{aH23} answered this question for cases when $s=n$, published posthumously in 1923. The same result was found independently by Radon \cite{jR22} in 1922.  Their result has been extended to any field $K$ in which $2 \ne 0$.

\HRThm \label{HRTheorem} An $[r, n, n]$-formula exists over $K$ if and only if $r \le \rho(n)$.

That maximal value $\rho(n)$, now called the Hurwitz-Radon function, is determined by the following rules:

\spb If $n = 1, 2, 4$ or 8 then $\rho(n) = n$.

\spb If $k$ is odd, then $\rho(2^mk) = \rho(2^m)$.

\spb $\rho(16n) = 8 + \rho(n)$.

It's easy to check that $\rho(n) = n$ only when $n = 1, 2, 4$ or 8, as expected from Hurwitz' earlier Theorem. Note that:\\[3pt]
\spb $ \rho(16) = 9, \quad \rho(32) = 10$, \quad and \quad  $\rho(64) = 12$.
   
The proof of Theorem \ref{HRTheorem} uses linear algebra, and we review the initial steps here.  Suppose an $[r, s, n]$-formula (*) is given as above.  Let $X$ be the column vector $(x_1, \dots, x_r)^\tr$, and similarly for column vectors $Y$ and $Z$. (Here $\tr$ denotes the transpose.)  Since $Z$ is bilinear in $X, Y$, we may  express it as $Z = (x_1A_1 + \cdots + x_rA_r)Y$ for some $n \times s$ matrices $A_i$ with entries in $K$.  With that notation, formula (*) is equivalent to the following system of ``Hurwitz Equations'' for those $n \times s$ matrixes $A_1, \dots, A_r$:

\pagebreak
\spd $A_i^\tr A_i = 1_s$\, \quad  whenever $1 \le i \le r$;  \\[5pt]
\spc $A_i^\tr A_j + A_j^\tr A_i = 0$ \; whenever $1 \le i,j \le r$, and $i \ne j$.

When $s = n$, those $A_j$ are square matrices and methods of linear algebra are easier to use. Further history and details about the Hurwitz-Radon Theorem appear in articles and books by several authors.  For instance, see the references in  \cite{dS00}.
 
Analysis of $[r, s, n]$-formulas is more difficult when $s < n$. In most cases there is a wide gap between sizes $(r, s, n)$ for which constructions are  known, and those that have been proved to be impossible.
  
The following ``Doubling Lemma''  is a first step in constructing such formulas.

\lemma \label{Double} An $[r, s, n]$-formula yields an $[r+1, 2s, 2n]$-formula.

\begin{proof} Here is a direct construction, as mentioned in \cite{dS00} Exercise 0.2.  From a given $[r, s, n]$ we have a system $A_1 \dots, A_r$ of $n \times s$ matrices satisfying the Hurwitz Equations.  Choose $A_1$ to get special treatment, and consider the following system of $(2n) \times (2s)$ matrices given in block form:

\spb $\left[\begin{smallmatrix} A_1 & 0 \\ 0 & \phantom{-}A_1 \end{smallmatrix}\right] , \spa \left[\begin{smallmatrix} \phantom{-}0 & A_1 \\- A_1 & 0 \end{smallmatrix}\right] , 
\spa  \left[\begin{smallmatrix} A_j & \phantom{-}0 \\ 0 & -A_j \end{smallmatrix}\right] $ \quad for $2 \le j \le r$ .  \\[5pt]
Some matrix multiplications verify that those $r+1$ matrices satisfy the Hurwitz Equations, yielding an $[r+1, 2s, 2n]$-formula.
\end{proof}

This proof involves the matrices 
$A_1 \otimes \left[\begin{smallmatrix} 1 & \phantom{-}0 \\ 0 & \phantom{-}1 \end{smallmatrix}\right] , \;
A_1 \otimes \left[\begin{smallmatrix}  \phantom{-}0 & \;1 \\ -1 & \;0 \end{smallmatrix}\right] $, \;and\;
$A_j \otimes \left[\begin{smallmatrix} 1 & \phantom{-}0 \\ 0 & -1 \end{smallmatrix}\right],$
where the symbol $\otimes$  stands for ``Kronecker product.''  For example, if $M$ is an $n \times s$ matrix, then $M \otimes  \left[\begin{smallmatrix}a & \;b \\  \rule{0mm}{2.5mm}c & \;d \end{smallmatrix}\right]  =  \left[\begin{smallmatrix} aM &\; bM \\ \rule{0mm}{2.5mm}cM & \;dM \end{smallmatrix}\right] $ is a $(2n) \times (2s)$ matrix given in block form.  On a more abstract level, this $\otimes$ arises from tensor products in the category of $K$-vector spaces.

Applying the Doubling Lemma three times to an $[8,8,8]$-formula, we obtain sizes:\\
\spb  $[9, 16, 16], \quad [10, 32, 32],$ \quad and \quad $[11, 64, 64]$.

Since $\rho(16) = 9$ and $\rho(32) = 10$ it does not seem that the $r+1$ entry in the Lemma can be improved.  But $\rho(64) = 12$, so there is a $[12, 64, 64]$, a larger size than the $[11, 64, 64]$ found by Doubling.  Such gaps can be filled in the classical case ($s = n$) by using the Expansion and Shift Lemmas described in \cite{dS00}. 
Recently,  C. Zhang and H.-L. Huang filled that gap a different way by improving the Doubling Lemma directly:

\proposition [Extended Doubling  \cite{ZH17}] \label{EDoub} \ \\
For any $k \ge 1$, an $[r, s, n]$-formula yields an $[r+\rho(2^{k-1}), \;2^ks, \;2^kn]$-formula.

Application of Zhang-Huang Doubling to the classical $[8,8,8]$-formula  yields a $[\rho(n), n, n]$-formula, for every $n = 2^m \ge 8$.

The goal of this note is to extend this Zhang-Huang result a step further, using the idea of ``amicable'' spaces.  As discussed in Chapter 2 of \cite{dS00}, that expanded view of compositions is useful in exposing some of the symmetries in the classical $[r, n, n]$-formulas.  Amicable spaces have not been investigated much when $s < n$. We begin with a definition: Matrices $A, B$ are \emph{amicable} if $A^\tr B = B^\tr A$.

Suppose a $[p, s, n]$ formula is given by the  $n \times s$ matrices $A_1, \dots, A_p$, and a $[q, s, n]$-formula is given by the $n \times s$ matrices $B_1, \dots, B_q$.  Then these systems satisfy the Hurwitz Equations

\spb $A_i^\tr A_i = 1_s$ and $A_i^\tr A_j + A_j^\tr A_i = 0$ \; for every $i \ne j$,\; and \\
\spb $B_k^\tr B_k = 1_s$ and $B_k^\tr B_\ell + B_\ell^\tr B_k = 0$ for every $k \ne \ell$.

Here we assume that $i, j$ run from 1 to $p$, while $k, \ell$ run from 1 to $q$.

\definition Those formulas of sizes $[p, s, n]$ and $[q, s, n]$ are \emph{amicable} if:  $A_i^\tr B_k = B_k^\tr A_i$ for every $i, k$.

An example of amicable formulas arose in the proof of Lemma \ref{Double}\,:  
Matrices
\;
 $\left[\begin{smallmatrix} 1 & \phantom{-}0 \\ 0 & \phantom{-}1 \end{smallmatrix}\right] , \; \left[\begin{smallmatrix}  \phantom{-}0 & \ 1 \\ -1 & \ 0 \end{smallmatrix}\right]$ , \quad and \quad
$\left[\begin{smallmatrix} 1 & \phantom{-}0 \\ 0 & -1 \end{smallmatrix}\right], \; \left[\begin{smallmatrix} 0 & \phantom{-}1 \\ 1 &\phantom{-} 0 \end{smallmatrix}\right] $ 
form a pair of amicable formulas each of size $[2, 2, 2]$.

The Doubling Lemma \ref{Double} extends to this context.

\lemma \label{AmicableDouble} Amicable $[p, s, n]$ and $[q, s, n]$ yield amicable  $[p+1, 2s, 2n]$ and $[q+1, 2s, 2n]$.

\begin{proof}
Start with matrices $A_i$ and $B_k$ as above.  
Consider the $2n \times 2s$ matrices:

\spb  $A_1 \otimes \left[\begin{smallmatrix} 1 & \phantom{-}0 \\ 0 & \phantom{-}1 \end{smallmatrix}\right] , \spz A_1 \otimes  \left[\begin{smallmatrix}  \phantom{-}0 & \ 1 \\ -1 & \ 0 \end{smallmatrix}\right] , \spz A_j \otimes \left[\begin{smallmatrix} 1 & \phantom{-}0 \\ 0 & -1 \end{smallmatrix}\right]$ \; for $2 \le j \le p$\,;\\[10pt]
\spb \spa $A_1 \otimes  \left[\begin{smallmatrix} 0 & \phantom{-}1 \\ 1 &\phantom{-} 0 \end{smallmatrix}\right]  , \spz B_k \otimes  \left[\begin{smallmatrix} 1 & \phantom{-}0 \\ 0 & -1 \end{smallmatrix}\right]$ \; for $1 \le k \le q$.

The $p+1$ matrices in the first row satisfy the Hurwitz Equations, and similarly for the $q+1$ matrices in the second row.  Moreover, the matrices in the first row are amicable with those in the second row.
\end{proof}

This doubling idea for amicable systems was pointed out in Exercise 2.12 of  \cite{dS00}. We apply it below only in the case $q=0$, when there are no matrices $B_k$. 

Here is the main result of this note.

 \theorem \label{Main} If   $[r, s, n]$ and $[r', s', n']$ formulas exist, then there is a $[r+r', \;2ss', \;2nn']$-formula.
 \medskip

Applying this when $[r', s', n'] = [\rho(2^{k-1}), 2^{k-1}, 2^{k-1}]$ yields Proposition \ref{EDoub}.

To begin, we investigate what conditions are needed to combine two formulas.  

Suppose an $[r, s, n]$-formula is given by the $n \times s$  matrices $A_1, \dots, A_r$, and $B$ is another matrix of that size. \\
Suppose a $[p, q, m]$-formula is given by the $m \times q$ matrices $C_1, \dots, C_p$, and $D$ is another matrix of that size. 

Consider the following system of $nm \times sq$ matrices:\\[5pt]
\spb  $A_j \otimes D$ for $1 \le j \le r$, \; and \; $B \otimes C_k$ for $1 \le k \le p$.  \\[5pt]
When does that list provide a $[r+p, \;sq, \;nm]$-formula?   

For $A_1 \otimes D, \;\dots\; , A_r \otimes D$ to satisfy the Hurwitz Equations, we need: $D^\tr D = 1_q$.  Similarly, for $B \otimes C_1, \;\dots\; , B \otimes C_p$ to satisfy those equations, we need:  $B^\tr B = 1_s$.  The remaining requirement is:\\[3pt]
\spb $(A_j^\tr B)  \otimes (D^\tr C_k) \;+\; (B^\tr A_j) \otimes (C_k^\tr D) \;=\; 0$ \; for every $j, k$.

Let's assume that $A_1, \dots, A_r, B$ form an $[r+1, s, n]$-formula.  That is:  $A_j^\tr B + B^\tr A_j = 0$ for every $j$. Then the condition above becomes: \\
\spc $(A_j^\tr B) \otimes \big(D^\tr C_k - C_k^\tr D\big) = 0$.  \\[5pt]
Then the Hurwitz Equations hold for the full list of $r+p$ matrices, provided:  $D$ is amicable with each $C_k$.

Here is a summary of what we have proved so far.

\proposition \label{FullDouble} Suppose amicable $[p, q, m]$ and $[1, q, m]$ formulas exist.  Then an $[r+1, s, n]$-formula yields an $[r+p, \;sq, \;nm]$-formula.

With those observations, the proof of  our Theorem is quickly done.

\begin{proof}[Proof of Theorem \ref{Main}\,.]    We are given a $[r', s', n']$ and an $[r, s, n]$. The Doubling Lemma \ref{AmicableDouble} applied to that $[r', s', n']$ produces amicable  $[r'+1, 2s', 2n']$ and $[1, 2s', 2n']$ formulas. Then Proposition \ref{FullDouble} yields the desired $[r+r', \;2ss', \;2nn']$-formula.
\end{proof}

\vsp

\bigskip

{\sc Daniel B. Shapiro, Department of Mathematics, Ohio State University} \\
Email: \href{shapiro@math.ohio-state.edu}{shapiro@math.ohio-state.edu}

\end{document}